\documentclass[a4paper,12pt]{article}
\usepackage[english]{babel}
\usepackage[latin1]{inputenc}
\usepackage{mathrsfs}
\usepackage{amsmath}
\usepackage{amsthm}
\usepackage{amssymb}
\usepackage{indentfirst}
\usepackage{type1cm}
\usepackage[dvips]{graphicx}

\usepackage{ifthen}
\newcommand{\cit}[1]{{\rm \textbf{#1}}}

\newcommand{\Ref}[2]{\cit{%
\ifthenelse{\equal{#1}{thm}}{Theorem}{}%
\ifthenelse{\equal{#1}{prop}}{Proposition}{}%
\ifthenelse{\equal{#1}{lem}}{Lemma}{}%
\ifthenelse{\equal{#1}{cor}}{Corollary}{}%
\ifthenelse{\equal{#1}{defn}}{Definition}{}%
\ifthenelse{\equal{#1}{oss}}{Remark}{}%
\ifthenelse{\equal{#1}{sec}}{Section}{}%
\ifthenelse{\equal{#1}{ex}}{Example}{}%
\ifthenelse{\equal{#1}{conj}}{Conjecture}{}%
\ifthenelse{\equal{#1}{ssec}}{Subsection}{}%
\ifthenelse{\equal{#1}{tab}}{Table}{}%
\ifthenelse{\equal{#1}{cla}}{Claim}{}%
\  \ref{#1:#2}%
}}

\theoremstyle{plain} 
\newtheorem{prop}{Proposition}[section]
\newtheorem{thm}[prop]{Theorem}
\newtheorem{lem}[prop]{Lemma} 
\newtheorem{cor}[prop]{Corollary}

\theoremstyle{remark}
\newtheorem*{proof1}{Proof of \Ref{prop}{pic20}}
\newtheorem*{proof2}{Proof of \Ref{thm}{sporadic}}
\newtheorem*{proof3}{Proof of \Ref{thm}{F_11}}

\newtheorem{oss}[prop]{Remark}
\newtheorem{ex}[prop]{Example}

\theoremstyle{definition}
\newtheorem{defn}[prop]{Definition}
\newcommand{\hk}{hyperk\"{a}hler }

\newcommand{\kahl}{K\"{a}hler }
\newcommand{\ktiposp}{$K3^{[2]}$ type }
\newcommand{\ktipo}{$K3^{[2]}$ type}
\title{On symplectic automorphisms of \hk fourfolds of \ktipo}
\author{Giovanni Mongardi}
\begin{document}
\maketitle
\abstract{The present paper proves that finite symplectic groups of automorphisms of \hk fourfolds deformation equivalent to the Hilbert scheme of two points on a $K3$ surface are contained in the simple group $Co_1$. Then we give an example of a symplectic automorphism of order 11 on the Fano scheme of lines of a cubic fourfold.}

\section{Introduction}
An automorphism $\varphi$ of a \hk manifold $X$ is symplectic if
\begin{equation}\nonumber
\varphi^*(\sigma_X)=\sigma_X,
\end{equation}
where $\sigma_X$ is a holomorphic symplectic 2-form on $X$.\\  
Finite abelian groups of symplectic automorphisms of complex K3 surfaces have been classified by Nikulin in \cite{nik1}. In particular one knows that a symplectic automorphism of finite order on a $K3$ surface over $\mathbb{C}$ has order at most 8.\\
The present paper deals with symplectic automorphisms on \hk fourfolds deformation equivalent to the Hilbert scheme of two points of a $K3$ surface. Such fourfolds will be called manifolds of \ktipo. 
Let us recall that manifolds of \ktiposp have $b_2=23$ and $H^2_{\mathbb{Z}}\cong U^3\,\oplus\,E_8(-1)^2\,\oplus\,(-2)$,  
where $U$ is the hyperbolic plane, $E_8(-1)$ is the unique negative definite even unimodular lattice of rank 8 and $(-2)$ is the rank 1 lattice of discriminant $-2$.\\
Let $Co_1$ be  Conway's sporadic simple group.
The main result of the paper is the following:
\begin{thm}\label{thm:sporadic}
Let $X$ be a \hk manifold of \ktiposp and let $G$ be a finite group of symplectic automorphisms of $X$. Then $G$ is isomorphic to a subgroup of $Co_1$. 
\end{thm}
We recall that Mukai \cite{muk} proved an analogous result for $K3$ surfaces (see also the proof of Kondo \cite{kon}): a finite group of symplectic automorphisms is a subgroup of Mathieu's group $M_{23}$.\\
A partial converse to \Ref{thm}{sporadic} is provided by \Ref{prop}{sfiga}, which also gives a computational method of determining possible finite automorphism groups.
Then we use \Ref{thm}{sporadic} to prove the following result on symplectic automorphisms of order 11:
\begin{prop}\label{prop:pic20}
Let $X$ be a fourfold of \ktiposp and let $\psi\,:\,X\,\rightarrow\,X$ be a symplectic automorphism of order 11. Then $21\,\geq\,h^{1,1}_{\mathbb{Z}}(X)\geq 20$. Moreover $\mathrm{Bir}(X)$ has a subgroup isomorphic to $\mathrm{PSL}_2(\mathbb{Z}_{/(11)})$. Furthermore if $X$ is projective then $h^{1,1}_{\mathbb{Z}}(X)=21$.
\end{prop}
Finally, we give an example of a fourfold of \ktiposp with a symplectic automorphism of order 11. Our example is given by the Fano scheme of lines on a cubic fourfold (unique up to projectivities). Fano schemes of lines on cubic fourfolds were first studied in \cite{bedon}, where the authors proved that they are of \ktipo.\\ 
The manifold $X_{Kl}\subset\mathbb{P}_\mathbb{C}^5$, given with respect to homogeneous coordinates $(x_0:\dots:x_5)$ by
\begin{equation}\label{fano_eq}
X_{Kl}\,=\,V(x_0^3+x_1^2x_5+x_2^2x_4+x_3^2x_2+x_4^2x_1+x_5^2x_3)
\end{equation}
is a $3:1$ cover of $\mathbb{P}_\mathbb{C}^4$, ramified along a cubic threefold first studied by Klein \cite{kle}. We denote the Fano scheme of lines in $X_{Kl}$ by $F_{Kl}$.
\begin{thm}\label{thm:F_11}
$F_{Kl}$ is a \hk fourfold of \ktipo. It has a symplectic automorphism $\varphi$ of order 11, induced by the element of $\mathrm{PGL}_6(\mathbb{C})$
\begin{equation}\label{phi_eq}
(x_0:x_1:x_2:x_3:x_4:x_5)\,\rightarrow\,(x_0:\omega x_1:\omega^3 x_2:\omega^4 x_3:\omega^5 x_4:\omega^9x_5),
\end{equation}
where $\omega=e^{\frac{2\pi i}{11}}$.
\end{thm}
We will describe other automorphisms of $X_{Kl}$. They are taken from the same work of Klein and were studied also by Adler \cite{adl}.\\
The paper is organized as follows. In \Ref{sec}{lattices} we recall some results in lattice theory and we use them to prove \Ref{thm}{sporadic} and \Ref{prop}{pic20}.\\
\Ref{sec}{polar_def} briefly analyses deformations of manifolds of \ktiposp having a symplectic automorphism of order 11. Moreover in it we also compute $\mathrm{NS}(X)$ for one interesting polarization.\\
In \Ref{sec}{ex_11} we give the promised example of a manifold with an order 11 symplectic automorphism, and in \Ref{sec}{5autom} we describe more symplectic automorphisms of this example.

\section{Lattice theory}\label{sec:lattices}
In this section we give a proof of \Ref{prop}{pic20} and a proof of \Ref{thm}{sporadic} using several results on general lattice theory and on lattices defined by symplectic automorphisms. The interested reader can consult \cite{nik2} for the main results concerning discriminant forms, \cite{con} for a broader treatment of \Ref{ssec}{niemeier} and \Ref{ssec}{holy}, and \cite{me1}
 for proofs of the stated results on lattices defined by symplectic automorphisms.\\

Let $L$ be a lattice, i.e. a free $\mathbb{Z}$-module equipped with an integer valued symmetric nondegenerate bilinear form $(\,,\,)_L$. We say that $L$ is even if $(a,a)_L\,\in\,2\,\mathbb{Z}$ for all $a\,\in\,L$, and unimodular if $L^\vee=L$.\\ Given an even lattice $L$ the group $A_L=L^{\vee}/L$ is called the discriminant group. Let $l(A_L)$ be the minimal number of generators of $A_L$. On $A_L$ there is a well-defined quadratic form $q_{A_L}$ (induced by $(\,,\,)_L$) taking values in $\mathbb{Q}/2\mathbb{Z}$ called the discriminant form. Let $(l_+,l_-)$ denote the signature of the quadratic form induced by $(\,,\,)_L$ on $L\otimes \mathbb{R}$.\\
With an abuse of notation the signature $\mathrm{sign}(q)$ of a discriminant form $q$ is defined as $l_+-l_-$ (modulo 8), where $L$ is a lattice with discriminant form $q$. This notion is well defined since two lattices $M,M'$ such that $q_{A_M}=q_{A_{M'}}$ are stably equivalent, i.e. there exist two unimodular lattices $T,T'$ such that $M\oplus T\cong M'\oplus T'$.\\
Two lattices $M$ and $M'$ are said to have the same genus if $M\otimes\mathbb{Z}_p\cong M'\otimes\mathbb{Z}_p$ for all primes $p$. Notice that there might be several isometry classes in the same genus.\\

Let $S$ be a nondegenerate sublattice of an unimodular lattice $L$ and let $M=S^\perp$. Then $A_M=A_S$ and $q_{A_M}=-q_{A_S}$.\\ 
If $G$ is a finite group of isometries of a lattice $L$, we let 
\begin{equation}\nonumber
T_G(L)=L^G
\end{equation}
be the invariant lattice and let 
\begin{equation}\nonumber
S_G(L)=T_G(L)^\perp
\end{equation}
be the co-invariant lattice.\\
The following are a simplified version of fundamental results on the existence of lattices and on the existence of primitive embeddings. Both general results were proven by Nikulin \cite[Theorem 1.10.1 and 1.12.2]{nik2}:
\begin{lem}\label{lem:nik_esiste}
Let $A_T$ be a finite abelian group and let $q_T$ be a quadratic form on $A_T$ with values in $\mathbb{Q}/2\mathbb{Z}$. Suppose the following are satisfied:
\begin{itemize}
\item there exists a lattice $T'$ of rank $t_++t_-$ and discriminant form $q_T$ over the group $A_T$;
\item $\mathrm{sign}(q_T)\equiv t_+-t_-$ $mod\,8$;
\item $t_+\geq0$, $t_-\geq0$ and $t_++t_-\geq l(A_T)$.
\end{itemize}
Then there exists an even lattice $T$ of signature $(t_+,t_-)$, discriminant group $A_T$ and form $q_{A_T}$.
\end{lem}
\begin{lem}\label{lem:nik_immerge1}
Let $S$ be an even lattice of signature $(s_+,s_-)$. There exists a primitive embedding of $S$ into some unimodular lattice $L$ of signature $(l_+,l_-)$ if and only if there exists a lattice $M$ of signature $(m_+,m_-)$ and discriminant form $q_{A_M}$ such that: 
\begin{itemize}
\item $s_++m_+=l_+$ and $s_-+m_-=l_-$;
\item $A_M\cong A_S$ and $q_{A_M}=-q_{A_S}$.
\end{itemize}
\end{lem}

\begin{oss}\label{oss:s_in_24}
Let $L=U^3\oplus E_8(-1)^2\oplus (-2)$ and let $G$ be a subgroup of $O(L)$. Then there exists a primitive embedding $L\rightarrow L'\cong U^4\oplus E_8(-1)^2$ such that $G$ extends to a group of isometries of $L'$ and $S_G(L)=S_G(L')$.
\begin{proof}
Let $x$ be a vector of square $2$ and $v\in L$ a vector of square $-2$ such that $(v,L)=2\mathbb{Z}$. Let $L'$ be the overlattice of $L\oplus\mathbb{Z}x$ generated by $L$ and $\frac{x+v}{2}$ and extend the action of $G$ to $L'$ by letting $G$ act as the identity on $x$. A direct computation shows $S_G(L)=S_G(L')$.
\end{proof}
\end{oss}
\subsection{Niemeier lattices and Leech couples}\label{ssec:niemeier}
In this subsection we recall Niemeier's list of negative definite even unimodular lattices in dimension 24 and we introduce a class of lattices which will be of fundamental interest in the rest of the section.
Detailed information about these lattices can be found in \cite[Chapter 16]{con} and in \cite[Section 1.14]{nik2}.\\
\begin{defn}
Let $M$ be a lattice and let $G\subset\,O(M)$. Then $(M,G)$ is Leech couple if the following are satisfied:
\begin{itemize}
\item $M$ is negative definite;
\item $M$ contains no vectors of square $-2$;
\item $G$ acts trivially on $A_M$;
\item $S_G(M)=M$.
\end{itemize} 
\end{defn}

Notice that $(\Lambda, Co_0)$ as in \Ref{tab}{nieme} is a Leech couple.\\
Now we recall Niemeier's list of negative definite even unimodular lattices of dimension 24. All of these lattices can be obtained by specifying a 0- or 24-dimensional negative definite Dynkin lattice such that every semisimple component has a fixed Coxeter number. In \Ref{tab}{nieme} we recall the possible choices. We obtain a Niemeier lattice $N$ from its Dynkin lattice $A$ by adding a certain set of glue vectors, which is a subset $G(N)$ of $A^{\vee}/A$. The precise definition of the glue vectors can be found in \cite[Section 4]{con} and we keep the notation from there. Notice that the set of glue vectors forms an additive subgroup of $A^{\vee}/A$.\\
The maximal Leech-type group $\mathrm{Leech}(N)$ is defined to be the maximal subgroup $G$ of $O(N)$ such that $(S_G(N),G)$ is a Leech couple. This group, first computed in \cite{ero}, is isomorphic to $O(N)/W(N)$, where $W(N)$ is the Weyl group generated by reflections in $-2$ vectors.\\
This data is summarized in \Ref{tab}{nieme}. In that table we have used the notation of \cite{atlas} for groups and the notation of \cite{con} for the glue code, which we recall for convenience:
\begin{itemize}
\item $n$ means a cyclic group of order $n$;
\item $p^n$ means an elementary $p$-group of order $p^n$;
\item $G.H$ means any group $F$ with a normal subgroup $G$ such that $F/G=H$;
\item $\mathrm{L}_m(n)$ means the group $\mathrm{PSL}_m$ over the finite field with $n$ elements;
\item $S_n$ and $A_n$ mean respectively permutation and alternating groups on $n$ elements;
\item $M_n$ and $Co_n$ are the Mathieu and Conway groups;
\end{itemize}
and, for the glue code,
\begin{itemize}
\item $[abc]$ means a vector $(x,y,z)$ with $x,y$ and $z$ glue vectors of type $a$, $b$ and $c$ respectively;
\item $[(abc)]$ means all glue vectors obtained from $[a,b,c]$ by cyclic permutations, hence $[abc],[bca],[cab]$.
\end{itemize}
It is a well known fact (see \cite[Chapter 26]{con}) that all the Niemeier lattices can be defined as sublattices of $\Pi_{1,25}\cong U\oplus E_8(-1)^3$ by specifying a primitive isotropic vector $v$ and setting $N=(v^\perp\cap\Pi_{1,25})/v$.

\begin{table}[ht]\label{tab:nieme}
\caption{Niemeier lattices}
\begin{tabular}{|c|c|c|c|c|}
 \hline
Name & Dynkin &Maximal Leech & Coxeter& Generating glue code\\ 
 & diagram & type Group &  Number & \\

\hline
$N_{1}$ & $D_{24}$ & $1$ & $46$ & $[1]$\\ \hline
$N_{2}$ &$D_{16}E_8$ & $1$ & $30$ & $[10]$\\ \hline
$N_{3}$ &$E_8^3$ & $S_3$ & $30$ & $[000]$\\ \hline
$N_{4}$ &$A_{24}$ & $2$ & $25$ & $[5]$\\ \hline
$N_{5}$ &$D_{12}^2$ & $2$ & $22$ & $[12],[21]$\\ \hline
$N_{6}$ &$A_{17}E_7$ & $2$ & $18$ & $[31]$\\ \hline
$N_{7}$ &$D_{10}E_7^2$ & $2$ & $18$ & $[110],[301]$\\ \hline
$N_{8}$ &$A_{15}D_9$ & $2$ & $16$ & $[21]$ \\ \hline
$N_{9}$ &$D_8^3$ & $S_3$ & $14$ & $[(122)]$ \\ \hline
$N_{10}$ &$A_{12}^2$ & $4$ & $13$ & $[15] $\\ \hline
$N_{11}$ &$A_{11}D_7E_6$ & $2$ & $12$ & $[111]$\\ \hline
$N_{12}$ &$E_6^4$ & $2.S_4$ & $12$ & $[1(012)]$ \\ \hline
$N_{13}$ &$A_9^2D_6$ & $2^2$  & $10$ & $[240],[501],[053]$ \\ \hline
$N_{14}$ &$D_6^4$ & $S_4$ & $10$ & $[$even perm. of $\{0,1,2,3\}]$ \\ \hline
$N_{15}$ &$A_8^3$ & $S_3\times 2$  & $9$ & $[(114)]$\\ \hline
$N_{16}$ &$A_7^2D_5^2$ & $2^3$ & $8$ & $[1112],[1721]$\\ \hline
$N_{17}$ &$A_6^4$ & $2.A_4$ & $7$ & $[1(216)]$ \\ \hline
$N_{18}$ &$A_5^4D_4$ & $2.S_4$ & $6$ & $[2(024)0],[33001],[30302],[30033] $ \\ \hline
$N_{19}$ &$D_4^6$ & $3\times S_6$ &$6$ & $[111111],[0(02332)]$ \\ \hline
$N_{20}$ &$A_4^6$ & $2.\mathrm{L}_2(5).2$ & $5$ & $[1(01441)]$\\ \hline
$N_{21}$ &$A_3^8$ & $2^3.\mathrm{L}_2(7).2$ & $4$ & $[3(2001011)]$ \\ \hline
$N_{22}$ &$A_2^{12}$ & $2.M_{12}$ &$3$ & $[2(11211122212)]$\\ \hline
$N_{23}$ &$A_1^{24}$ &  $M_{24}$ & $2$ & $[1(00000101001100110101111)]$ \\ \hline
$\Lambda$ &$\emptyset$ & $Co_0$ &$0$ & $\emptyset$ \\ \hline
\end{tabular}
\end{table}

\begin{ex}
Let $\Pi_{1,25}\subset\mathbb{R}^{26}$ (the first coordinate of $\mathbb{R}^{26}$ is the positive definite one) be as before and let
\begin{eqnarray}\nonumber
v &= &(17,1,1,1,1,1,1,1,1,3,3,3,3,3,3,3,3,3,5,5,5,5,5,5,5,5)\\\nonumber
w &= &(70,0,1,2,3,4,5,\dots,24)
\end{eqnarray}
be two isotropic vectors in the standard basis of $\mathbb{R}^{26}$. Then 
\begin{equation}\nonumber
\Lambda\cong (w^\perp\cap\Pi_{1,25})/w
\end{equation}
 and 
\begin{equation}\nonumber 
N_{15}\cong (v^\perp\cap\Pi_{1,25})/v.
\end{equation}
\end{ex}

\subsection{The ``holy" construction and\\ automorphisms of the Leech lattice}\label{ssec:holy}
In this subsection we sketch the ``holy construction" (see \cite[Chapter 24]{con} for details) of the Leech lattice $\Lambda$ from other Niemeier lattices. We shall use it later in proving \Ref{prop}{pic20}.
\\
$A_n(-1)$ is the negative definite Dynkin lattice defined by
\begin{eqnarray}\nonumber
A_n&=&\{(a_1,\dots,a_{n+1})\,\in\,\mathbb{Z}^{n+1},\,\,\sum{a_i}=0\},\\\nonumber
q_{A_n(-1)}&=&-q_{A_n}.
\end{eqnarray}
Let $f_0$ be the vector with $-1$ in the first coordinate and $1$ in the second, zero otherwise. Let $g_0=h^{-1}(-\frac{1}{2}n,-\frac{1}{2}n+1,\dots,\frac{1}{2}n)$ where $h$ is the Coxeter number of $A_n$. Let $f_j$ and $g_j,\,j\in\,\{1,\dots n\}$ be the images of $f_0$ and $g_0$ respectively under cyclic permutations of coordinates. Notice that $f_j,\,j\neq0$ are the simple roots of $A_n(-1)$.\\
Suppose now that $nm=24$, so that $A_n(-1)^m$ is a 24-dimensional lattice contained in a Niemeier lattice $N$, and let $h_k=(g_{j_1},\dots,g_{j_m})$, where $[j_1j_2\dots j_m]$ is a glue code obtained from \Ref{tab}{nieme} and $k\in G(N)$. Let $f^r_j=(0,\dots,0,f_j,0,\dots,0)$ where $f_j$ belongs to the $r$-th copy of $A_n(-1)$. Let $m_j^r$ and $n_k$ be integers.\\
\begin{prop}\cite[Chapter 24]{con}
With notation as above
\begin{equation}\label{holy_nieme}
\{\sum_{r=1}^m\,\sum_{j=1}^n m_j^rf_j^r+\sum_{k\in G(N)} n_kh_k,\,\,\text{such that }\,\,\,\,\sum_{k}n_k=0\} 
\end{equation}
is isometric to the negative definite Niemeier lattice $N$. 
\begin{equation}\label{holy_leech}
\{\sum_{r=1}^m\,\sum_{j=1}^n m_j^rf_j^r+\sum_{k\in G(N)} n_kh_k,\,\text{such that }\,\,\,\,\sum_{k}n_k+\sum_{j,r}m_j^r=0\}
\end{equation}
is isometric to the Leech lattice $\Lambda$ and is called the ``holy" construction of $\Lambda$ with hole $N$.\\
\end{prop}
Moreover the glue code provides several automorphisms of the Leech lattice, where the action of $t\in G(N)$ is given by sending $h_w$ to $h_{w+t}$.\\
The holy construction can be used to exhibit the action of certain elements of $Co_1$ on $\Lambda$ explicitly, as in the following examples.
\begin{ex}\label{ex:A122}
Let us apply this construction to the lattice $A_{12}^{2}$, where $G(N)=\mathbb{Z}_{/(13)}$. Let $\chi$ be an  automorphism of $\Lambda$ of order 13 generated by a nontrivial element $g$ of $G(N)$ via this holy construction. $\chi$ cyclically permutes the simple roots of both copies of $A_{12}$ and therefore has no fixed points in $\Lambda$.
\end{ex}
\begin{ex}\label{ex:p11A212}
Consider the lattice $A_2^{12}$ and the cyclic permutation $\chi$ of the last 11 copies. This defines automorphisms, also denoted $\chi$, of both $N_{22}$ and $\Lambda$ via the same action on glue vectors.\\
A direct computation shows $T_{\chi}(N_{22})$ is spanned by
\begin{equation}\nonumber
f_1^1,\,\,f_1^2,\,\,\sum_2^{12}f_1^i,\,\,\sum_2^{12}f_2^i,\,\,\sum_1^{11} h_j,
\end{equation}
where $h_j, j\in\{1,\dots,11\}$ are obtained from generators of the glue code as in \Ref{tab}{nieme}. Moreover $S_{\chi}(N_{22})$ has rank 20 and is spanned by
\begin{equation}
(f^k_1-\chi f^k_1),(f^k_2-\chi f^k_2),(g_j-\chi g_j).
\end{equation}
Here $k$ runs from 2 to 12 and $j$ is as before. These vectors lie in the set defined by  \eqref{holy_leech}, where $g_j$ plays the role of $h_k,k\in\,G(N)$. Therefore $S_{\chi}(N_{22})$ is contained in $\Lambda$ and, since both $S_{\chi}(N_{22})$ and $S_{\chi}(\Lambda)$ are primitive, $S_{\chi}(N_{22})=S_{\chi}(\Lambda)$.
\end{ex}
\begin{ex}\label{ex:p11A124}
A similar computation can be done for $A_1^{24}$. We index the copies of $A_1$ by the set
\begin{equation}\nonumber
\{\infty,0,1,\dots,22\}=\mathbb{P}^1(\mathbb{Z}_{/(23)}).
\end{equation}
The isometry $\chi$ of order 11 is defined by the permutation:
\begin{equation}\label{p11}
(0)(15\,7\,14\,5\,10\,20\,17\,11\,22\,21\,19)(\infty)(3\,6\,12\,1\,2\,4\,8\,16\,9\,18\,13).
\end{equation}
As before this isometry preserves both $N_{23}$ and $\Lambda$, and the lattice $S_{\chi}(N_{23})$ is generated by the following vectors:
\begin{equation}
(f^k_1-\chi f^k_1),(f^l_1-\chi f^l_1),(g_j-\chi g_j).
\end{equation}
Here $k$ runs through the indexes contained in the first 11-cycle of \eqref{p11}, $l$ runs  through the second one and $j$ through the generators of the glue code contained in \Ref{tab}{nieme}.\\
Once again all of these generators also lie in $\Lambda$, so
$S_{\chi}(N_{23})=S_{\chi}(\Lambda)$.
A direct computation shows that the lattice $S_{11}=S_{\chi}(N_{23})$ is given by the following quadratic form:

\begin{equation}\nonumber
\tiny{
\left(
\begin{array}{cccccccccccccccccccc}
 -4 & 1 & -2 & -2 & -1 & 1 & -1 & 1 & -1 & -1 & 2 & 1 & -1 & 2 & -1 & -2 & -2 & 2 & 1 & -1 \\
 1 & -4 & -1 & -1 & -1 & -1 & -1 & 1 & -1 & 2 & -1 & -2 & 2 & 0 & -1 & 0 & 0 & -1 & -2 & 1 \\
 -2 & -1 & -4 & -2 & -1 & -1 & 0 & 1 & 0 & -1 & 1 & 0 & -1 & 2 & -2 & -1 & -1 & 0 & 0 & 1 \\
 -2 & -1 & -2 & -4 & 0 & 0 & -2 & 0 & -1 & 0 & 2 & 1 & 0 & 1 & 0 & 0 & -1 & 1 & 0 & -1 \\
 -1 & -1 & -1 & 0 & -4 & 1 & -1 & 2 & -2 & -1 & 1 & 0 & -1 & 0 & -2 & -2 & 0 & 1 & 1 & -1 \\
 1 & -1 & -1 & 0 & 1 & -4 & 0 & -1 & 0 & 1 & -2 & -1 & 0 & -1 & -1 & 0 & -1 & 0 & -1 & 1 \\
 -1 & -1 & 0 & -2 & -1 & 0 & -4 & 1 & -2 & 1 & 1 & 1 & 0 & -1 & 0 & -1 & 0 & 2 & 0 & -2 \\
 1 & 1 & 1 & 0 & 2 & -1 & 1 & -4 & 0 & 0 & -1 & 1 & 1 & 0 & 2 & 1 & 0 & -1 & 1 & 0 \\
 -1 & -1 & 0 & -1 & -2 & 0 & -2 & 0 & -4 & 0 & 0 & 1 & 1 & 0 & -1 & -2 & 0 & 2 & 0 & -2 \\
 -1 & 2 & -1 & 0 & -1 & 1 & 1 & 0 & 0 & -4 & 1 & 1 & -2 & 1 & 0 & 0 & 1 & 1 & 1 & 0 \\
 2 & -1 & 1 & 2 & 1 & -2 & 1 & -1 & 0 & 1 & -4 & -2 & 2 & -1 & 0 & 0 & 0 & -1 & -2 & 1 \\
 1 & -2 & 0 & 1 & 0 & -1 & 1 & 1 & 1 & 1 & -2 & -4 & 1 & 0 & -1 & 0 & -1 & -1 & -2 & 2 \\
 -1 & 2 & -1 & 0 & -1 & 0 & 0 & 1 & 1 & -2 & 2 & 1 & -4 & 0 & -1 & 0 & 0 & 1 & 2 & 0 \\
 2 & 0 & 2 & 1 & 0 & -1 & -1 & 0 & 0 & 1 & -1 & 0 & 0 & -4 & 1 & 1 & 1 & 0 & 0 & -1 \\
 -1 & -1 & -2 & 0 & -2 & -1 & 0 & 2 & -1 & 0 & 0 & -1 & -1 & 1 & -4 & -2 & -1 & 1 & 0 & 0 \\
 -2 & 0 & -1 & 0 & -2 & 0 & -1 & 1 & -2 & 0 & 0 & 0 & 0 & 1 & -2 & -4 & -2 & 2 & 0 & -1 \\
 -2 & 0 & -1 & -1 & 0 & -1 & 0 & 0 & 0 & 1 & 0 & -1 & 0 & 1 & -1 & -2 & -4 & 1 & 0 & 0 \\
 2 & -1 & 0 & 1 & 1 & 0 & 2 & -1 & 2 & 1 & -1 & -1 & 1 & 0 & 1 & 2 & 1 & -4 & 0 & 2 \\
 1 & -2 & 0 & 0 & 1 & -1 & 0 & 1 & 0 & 1 & -2 & -2 & 2 & 0 & 0 & 0 & 0 & 0 & -4 & 1 \\
 -1 & 1 & 1 & -1 & -1 & 1 & -2 & 0 & -2 & 0 & 1 & 2 & 0 & -1 & 0 & -1 & 0 & 2 & 1 & -4
\end{array}
\right)
}
\end{equation}
It is worth mentioning that $S_{11}$ is in the same genus as $E_8(-1)^2\,\oplus\,\left(\begin{array}{cc} -2 & 1\\ 1 & -6 \end{array}\right)^2$ and $D_{16}^+(-1)\oplus\left(\begin{array}{cc} -2 & 1\\ 1 & -6 \end{array}\right)^2$, where $D_{16}^+$ is a unimodular overlattice of the Dynkin lattice $D_{16}$. 

\end{ex}
\subsection{Finite symplectic automorphism groups}
Let $X$ be a \hk manifold and let $G\subset \mathrm{Aut}(X)$. Then we put $S_G(X)=S_G(H^2(X,\mathbb{Z}))$ and $T_G(X)=T_G(H^2(X,\mathbb{Z}))$.
The following lemmas are contained in \cite{me1}:
\begin{lem}\label{lem:algaction}
If $G$ is a finite group of symplectic automorphisms of a fourfold $X$ of \ktipo, then
\begin{enumerate}
\item $S_G(X)$ is nondegenerate and negative definite;
\item $S_G(X)$ contains no elements with square $-2$;
\item $S_G(X)\subset Pic(X)$;
\item $G$ acts trivially on $A_{S_G(X)}$.
\end{enumerate}
\end{lem}
Notice that this amounts to saying that $(S_G(X),G)$ is a Leech couple. 
\begin{lem}\label{lem:cohom_to_bir}
Let $L=U^3\oplus E_8(-1)^2\oplus (-2)$ and let $G$ be a finite subgroup of $O(L)$. Suppose the following hold:
\begin{enumerate}
\item $S_G(L)$ is nondegenerate and negative definite;
\item $S_G(L)$ contains no element with square $(-2)$.
\end{enumerate}
Then there exists a \hk manifold $X$ of \ktiposp and a subgroup $G'\subset \mathrm{Bir}(X)$ such that $G'\cong G$, $S_G(L)\cong S_{G'}(X)$ and\\ $G'_{|H^{2,0}(X)}=Id$.
\end{lem}

\begin{prop}\label{prop:sfiga}
Let $(S,G)$ be a Leech couple such that $S$ is primitively contained in some Niemeier lattice $N$, and suppose there exists a primitive embedding $S\rightarrow L$.\\
Then $G$ extends to a group of bimeromorphisms on some \hk manifold $X$ of \ktipo.
\begin{proof}
This is an immediate consequence of \Ref{lem}{cohom_to_bir}: $G$ acts trivially on $A_S$, therefore we can extend $G$ to a group of isometries of $L$ acting trivially on $S^{\perp_L}$. Thus we have $S_G(L)\cong S$. Then the conditions of \Ref{lem}{cohom_to_bir} are satisfied because $(S,G)$ is a Leech couple.
\end{proof}
\end{prop}
We are now ready to prove the main result of this section:

\begin{proof2}
Let $b=\mathrm{rank}(S_G(X))$; by \Ref{lem}{algaction} $S_G(X)$ has signature $(0,b)$. By \Ref{oss}{s_in_24} we have a lattice $T'$ of signature $(4,20-b)$ such that $A_{T'}=A_{S_G(X)}$ and $q_{T'}=-q_{A_{S_G(X)}}$. Therefore we can apply \Ref{lem}{nik_esiste} obtaining a lattice $T$ of signature $(0,24-b)$ and discriminant form $-q_{A_{S_G(X)}}$. Thus by \Ref{lem}{nik_immerge1} there exists a primitive embedding $S_G(X)\rightarrow N$, where $N$ is one of the lattices contained in \Ref{tab}{nieme}. Again by \Ref{lem}{algaction} we see that $(S_G(X),G)$ is a Leech couple, hence $G$ lies inside $\mathrm{Leech}(N)$. Using the holy construction we obtain $\mathrm{Leech}(N)\subset \mathrm{Leech}(\Lambda)$. Moreover a direct computation shows that for all $G\subset \mathrm{Leech}(N)\subset \mathrm{Leech}(\Lambda)$ we have $\mathrm{rank}(S_G(N))=\mathrm{rank}(S_G(\Lambda))$ (after tensoring with $\mathbb{Q}$ they are both generated by elements of the form $v-g(v)$, $v\in N$ and $g\in G$). Obviously the central involution of $Co_0$ has a co-invariant lattice of rank 24, hence we can restrict ourselves to $Co_1$. 
\end{proof2}
\begin{cor}\label{cor:max_p_order}
Let $\phi$ be a symplectic automorphism of prime order $p$ on a \hk fourfold $X$ of \ktipo. Then $p\leq 11$.
\begin{proof}
By \Ref{thm}{sporadic} the order of a symplectic automorphism must divide the order of the group $Co_1$. That excludes all primes apart from $2,3,5,7,11,13,23$. An automorphism of order 23 has a co-invariant lattice which is negative definite and of rank 22, and therefore cannot embed into $H^2(X,\mathbb{Z})$. This can be explicitly computed using an order 23 element of $M_{24}$ and letting it act on $\Lambda$ or on $N_{23}$. The only Niemeier lattice with an automorphism of order 13 is $\Lambda$, where all elements of order 13 are conjugate (see \cite{atlas}). These automorphisms have no fixed points on $\Lambda$, as in \Ref{ex}{A122}.

\end{proof}
\end{cor}
\begin{proof1}
By \Ref{thm}{sporadic} $S_\psi(X)$ embeds in a Niemeier lattice $N$ and $\psi$ extends to an element of $\mathrm{Leech}(N)$. By \Ref{tab}{nieme} $N$ can only be $N_{22},N_{23}$ or $\Lambda$ and, up to conjugacy, there is only one possible choice for $\psi\in O(N_{23})$ or $\psi\in O(\Lambda)$ and there are only two possible choices for $\psi\in O(N_{22})$.\\ However we computed $S_{\psi}(N)$ in \Ref{ex}{p11A212} and \Ref{ex}{p11A124}, where we proved that it is isometric to $S_{11}$. We immediately have $20\leq h^{1,1}_{\mathbb{Z}}(X)\leq 21$ and obviously if $X$ is projective $h^{1,1}_{\mathbb{Z}}(X)=21$. Now we wish to give an action of $\mathrm{L}_2(11)=\mathrm{PSL}_2(\mathbb{Z}_{/(11)})$ on $S_{11}$. Without loss of generality we can suppose $N=N_{23}$ and the action of $\mathrm{L}_2(11)$ is given by the following permutations of the standard coordinates of $N_{23}$ (see \cite[pages 274 and 280]{con}):
\begin{eqnarray}\nonumber
\alpha &= &(15\,7\,14\,5\,10\,20\,17\,11\,22\,21\,19)(3\,6\,12\,1\,2\,4\,8\,16\,9\,18\,13),\\
\nonumber\beta&=&(14\,17\,11\,19\,22)(20\,10\,7\,5\,21)(18\,4\,2\,6\,1)(8\,16\,13\,9\,12),\\
\nonumber\gamma&=&(2\,4)(5\,10)(6\,18)(8\,12)(9\,16)(11\,17)(14\,19)(20\,21).
\end{eqnarray} 
Now by \Ref{prop}{sfiga} this action of $\mathrm{L}_2(11)$ on $S_{11}$ is induced by a group of birational transformations isomorphic to $\mathrm{L}_2(11)$.
\end{proof1}
\section{Deformation behaviour}\label{sec:polar_def}

In this section we analyse deformation classes of manifolds of \ktiposp with a symplectic automorphism of order 11 and we look at their possible invariant polarizations. 
\begin{defn}
Let $X$ be a \hk manifold with \kahl class $\omega$ and symplectic form $\sigma_X$. Then there exists a family
\begin{eqnarray} 
TW_{\omega}(X) &:= &X\times\mathbb{P}^1\\\nonumber
& & \downarrow\\\nonumber
S^2&\cong&\mathbb{P}^1 
\end{eqnarray}
called the twistor family, such that $TW_{\omega}(X)_{(a,b,c)}=X$ with complex structure given by the \kahl class $a\omega+b(\sigma_X+\overline{\sigma}_X)+c(\sigma_X-\overline{\sigma}_X)$.
\end{defn}
\begin{oss}\label{oss:twistor_autom}
Let $X$ be a \hk manifold and let $G\subset \mathrm{Aut}(X)$ be a finite group of symplectic automorphisms. Let $\omega$ be a $G-$invariant \kahl class. Then the action of $G$ on $X$ extends to a symplectic action of $G$ on all the fibers of the twistor space associated to $\omega$.
\begin{proof}
Every fiber has a \kahl class which is a linear combination of $\sigma_X,\overline{\sigma_X}$ and $\omega$. Since $G$ is symplectic on $X$ these classes are all $G-$invariant, hence $G\subset \mathrm{Aut}(TW_{\omega}(X)_t)$ for all $t$. 
\end{proof}
\end{oss}
Let $X$ be a manifold of \ktiposp with a symplectic automorphism $\psi$ of order 11 and let $\omega$ be a $\psi$-invariant \kahl class.\\ First, notice that it follows from \Ref{prop}{pic20} that a nontrivial deformation of $(X,\psi)$ has dimension at most 1. Moreover the twistor family $TW_{\omega}(X)$ is naturally endowed with a symplectic automorphism of order 11 as in \Ref{oss}{twistor_autom}.   
Therefore $TW_{\omega}(X)$ is already a family of the maximal dimension for such pairs $(X,\psi)$. Moreover we have that the twistor family $TW_{\omega}(X)$ is actually a family over the base $(T_{\psi}(X)\otimes\mathbb{R})/\mathbb{C}^*$, where the $\mathbb{C}^*$ action is given by the identification $T_{\psi}(X)=\langle\omega,\sigma_X,\overline{\sigma}_X\rangle\cap H^2_{\mathbb{Z}}(X)$.\\ Thus what we really need to analyse are the possible lattices $T_{\psi}(X)$ up to isometry. We have already proved that there exists only one isometry class of lattices $S_{\psi}(X)$, namely that of $S_{11}$. However there might be several isometry classes of lattices $T_{\psi}(X)$. In fact \Ref{thm}{sporadic} and \Ref{prop}{sfiga} can be used only to compute the genus of $T_{\psi}(X)$.\\ A direct computation shows that there are two such lattices, namely:
\begin{equation}\label{TX_1}
T^1_{11}=\left(\begin{array}{ccc}  2&1&0\\ 1&6&0\\ 0&0&22\end{array}\right),
\end{equation}
\begin{equation}\label{TX_2}
T^2_{11}=\left(\begin{array}{ccc} 6&-2&-2\\ -2& 8&-3\\-2&-3&8\end{array}\right).
\end{equation}
Therefore there are two distinct families of \hk manifolds endowed with a symplectic automorphism of order 11, whose existence is a consequence of \Ref{lem}{cohom_to_bir}. We call them $TW(X_1)$ and $TW(X_2)$. 
\subsection{Invariant Polarizations}
In this subsection we look at possible invariant polarizations of small degree of $TW(X_1)$ and $TW(X_2)$, i.e. at primitive vectors in $T^1_{11}$ and $T^2_{11}$. We computed several polarizations in degree up to 24.
\begin{prop}\label{prop:polar_TW1}
The minimal degree of an invariant polarization inside $TW(X_1)$ is $2$, and there are no polarizations in degrees $4,12,14,16$ or $20$. Moreover the least degree of a polarization $f$ such that $(f,L)=2\mathbb{Z}$ (i.e. $f$ has divisor 2) is $22$. 
\begin{proof}
Let $a,b,c$ be the basis of $T^1_{11}$ in \eqref{TX_1}. A minimal polarization is given by the vector $a$; a minimal polarization of degree $22$ and divisor 2 is given by $c$. The rest is just a direct computation.
\end{proof}
\end{prop}
\begin{prop}\label{prop:polar_TW2}
The minimal degree of an invariant polarization inside $TW(X_2)$ is $6$, and there are no polarizations in degrees $12,14,16$ or $20$. Moreover the least degree of a polarization $f$ such that $(f,L)=2\mathbb{Z}$ is $6$. 
\begin{proof}
Let $a,b,c$ be the basis of $T^2_{11}$ in \eqref{TX_2}. A minimal polarization is given by the vector $a$ which has also divisor 2. 
 The rest is just a direct computation.
\end{proof}
\end{prop}
\begin{cor}\label{cor:ns_11}
Let $X$ be a manifold of \ktiposp with a symplectic automorphism $\psi$ of order 11 and an invariant polarization of degree 6 and divisor 2. Then 
\begin{equation}
NS(X)\cong (6)\,\oplus\, E_8(-1)^2\,\oplus\,\left(\begin{array}{cc} -2 & 1\\ 1 & -6 \end{array}\right)^2
\end{equation}
and
\begin{equation}
T(X)=\langle\sigma_{X},\overline{\sigma_{X}}\rangle\,\cap\, H^2(X,\mathbb{Z})\cong\left(\begin{array}{cc} 22 & 33\\ 33&66 \end{array}\right).
\end{equation}
\begin{proof}
$NS(X)$ is an overlattice of $S_{11}\oplus (6)$. Since there are no nontrivial overlattices (its discriminant group has no nontrivial isotropic elements) of $S_{11}\oplus (6)$ we have $NS(X)=S_{11}\oplus (6)$. A direct computation shows that this lattice is isomorphic to $(6)\,\oplus\, E_8(-1)^2\,\oplus\,\left(\begin{array}{cc} -2 & 1\\ 1 & -6 \end{array}\right)^2$.\\
Finally, $T(X)$ is the orthogonal complement of the polarization in $T_{\psi}(X)$. By \Ref{prop}{polar_TW1} and \Ref{prop}{polar_TW2} we have $T_{\psi}(X)=T^2_{11}$ and a direct computation shows $T(X)=\left(\begin{array}{cc} 22 & 33\\ 33&66 \end{array}\right)$
\end{proof}
\end{cor}


\section{The Fano scheme of lines $F_{Kl}$}\label{sec:ex_11}
The main goal of this section is to prove \Ref{thm}{F_11}. We start with some results about Fano schemes of lines on cubic fourfolds, due to Beauville and Donagi \cite{bedon}.\\
\begin{thm}\label{thm:fano_map}
Let $X\,\subset\,\mathbb{P}^5$ be a smooth cubic fourfold and let $F(X)$ be the scheme parametrizing lines contained in $X$. Then the following hold:
\begin{itemize}
\item $F(X)$ is a \hk manifold.
\item $F(X)$ is of \ktipo.
\item the Abel-Jacobi map
\begin{equation}\label{abel} \alpha\,:\,H^4(X,\mathbb{C})\,\rightarrow\,H^2(F(X),\mathbb{C})
\end{equation}
is an isomorphism of Hodge structures.
\end{itemize}
\end{thm} 

\begin{proof3}
Let $\psi$ be the element of $\mathrm{PGL}_6(\mathbb{C})$ sending $(x_0:x_1:x_2:x_3:x_4:x_5)$ to $(x_0:\omega x_1:\omega^3 x_2:\omega^4 x_3:\omega^5 x_4:\omega^9x_5)$, where $\omega=e^{\frac{2\pi i}{11}}$. A cubic polynomial is $\psi$-invariant if and only if it is in the linear span of 
\begin{equation}\nonumber B=\{x_0^3,x_1^2x_5,x_2^2x_4,x_3^2x_2,x_4^2x_1,x_5^2x_3\},
\end{equation}
An easy computation shows that the differential of $f$ has nontrivial zeroes (so that the cubic fourfold $V(f)$ is singular) if and only if $f$ lies in the span of some proper subset of $B$. This is not the case for
\begin{equation}\nonumber 
h=x_0^3+x_1^2x_5+x_2^2x_4+x_3^2x_2+x_4^2x_1+x_5^2x_3,
\end{equation}
Therefore $X_{kl}=V(h)$, as defined in \eqref{fano_eq}, is smooth, and we can apply \Ref{thm}{fano_map}. We obtain that $F_{Kl}$ is a \hk manifold deformation equivalent to $K3^{[2]}$. Moreover we have the Hodge isomorphism $\alpha$ given in \eqref{abel}.\\
Let $\varphi$ be the map induced on $F_{Kl}$ by $\psi$. Using $\alpha$ we have that $\varphi$ is symplectic if and only if $\psi_{|H^{3,1}(X_{Kl})}=Id$.\\
By the formula in \cite[Théorème 18.1]{voi} we have
\begin{equation}\label{sympl}
H^{3,1}(X_{Kl})\,=\,\langle\mathrm{Res}(\frac{\Omega}{h^2})\rangle,
\end{equation}
where $\Omega=\sum_{i}\,(-1)^ix_i\,dx_0\wedge\,\dots\,\widehat{dx_i}\,\dots\,\wedge dx_5$. Since $\psi$ acts trivially on both $\Omega$ and $h$ we obtain that $\varphi(\sigma_{F_{Kl}})=\sigma_{F_{Kl}}$ for any symplectic 2-form $\sigma_{F_{Kl}}$ on $F_{Kl}$.
Finally fixed points of $\psi$ on $X_{Kl}$ are the eigenvectors of $\psi$ lying on $X_{Kl}$, which are the points $[e_1],[e_2],[e_3],[e_4],[e_5]$, where $e_1=(0,1,0,0,0,0)\in\mathbb{C}^6$, etc. A fixed line on $X_{Kl}$ must contain two fixed points, and therefore the fixed points on $F_{Kl}$ are those parametrizing lines through those points, namely the five lines
\begin{equation}\nonumber
\overline{[e_1][e_2]},\overline{[e_1][e_3]},\overline{[e_2][e_5]},\overline{[e_3][e_4]},\overline{[e_4][e_5]}.
\end{equation}
This proves that $\varphi$ is not the identity, hence it has indeed order 11.
\end{proof3}
\begin{oss}
\Ref{prop}{polar_TW1} and \Ref{prop}{polar_TW2} imply that $F_{Kl}\subset TW(X_2)$, and  \Ref{cor}{ns_11} gives the Neron-Severi and transcendental lattices. Moreover $\varphi$ does not lift to any small projective non-trivial deformation of $F_{Kl}$, since $h^{1,1}_{\mathbb{Z}}(F_{Kl})=h^{1,1}_{\mathbb{C}}(F_{Kl})$.

\end{oss}
\section{$\mathrm{L}_2(11)$ acting on $F_{Kl}$}\label{sec:5autom}
In this subsection we exhibit directly the automorphisms of $X_{Kl}$ that generate the action of $\mathrm{L}_2(11)$ from \Ref{prop}{pic20}.\\ Let us remark that $X_{Kl}$ is a 3 to 1 Galois cover of $\mathbb{P}^4$ ramified along the threefold 
\begin{equation}
KA=V(x_1^2x_5+x_2^2x_4+x_3^2x_2+x_4^2x_1+x_5^2x_3),
\end{equation}
where the covering map is simply the projection 
\begin{equation}\nonumber
[x_0,x_1,x_2,x_3,x_4,x_5]\,\rightarrow\,[x_1,x_2,x_3,x_4,x_5],
\end{equation}
 and the covering automorphism group is generated by 
 \begin{equation}\nonumber
 [x_0,x_1,x_2,x_3,x_4,x_5]\,\stackrel{\alpha}{\rightarrow}\,[\eta x_0,x_1,x_2,x_3,x_4,x_5],\,\,\,\eta=e^{\frac{2\pi i}{3}}.
 \end{equation} 
Notice that $\alpha$ acts as multiplication by $\eta$ on $H^{3,1}(X_{Kl})$.\\
Obviously any automorphism of $\mathbb{P}^4$ that preserves $KA$ extends to an automorphism of $X_{Kl}$. By the results in \cite{adl} and \cite{kle} these automorphism generate precisely the group $\mathrm{L}_2(11)\,=\,\mathrm{PSL}_2(\mathbb{Z}_{/(11)})$, which is a finite simple group of order 660.
Hence the automorphism group of $X_{Kl}$ contains $\mathbb{Z}_{/(3)}\,\times\,\mathrm{L}_2(11)$. Now we only need to find generators of this group and to determine whether they act symplectically on $F_{Kl}$ or not.\\
Permuting the coordinates on $\mathbb{P}^4$ by the cyclic permutation $(1\,4\,2\,3\,5)$ preserves $KA$ and therefore induces an automorphism $\beta$ of order 5 on $F_{Kl}$.
From \eqref{sympl} we can check that $\beta$ is a symplectic automorphism. Furthermore a direct computation on the Jacobian ring $\mathbb{C}[x_0,x_1,x_2,x_3,x_4,x_5]/(\partial h/\partial x_0,\dots,\partial h/\partial x_5)$ of $X_{Kl}$ shows that\\ $\mathrm{rank}(S_\beta(F_{Kl}))=16$.\\
Let $H$ be the kernel of the action of $\mathrm{L}_2(11)$ on $H^{2,0}(F_{Kl})$. $H$ is nontrivial since $\beta\in H$, therefore $H=\mathrm{L}_2(11)$. Hence $\mathrm{L}_2(11)$ acts symplectically on $F_{Kl}$.

$\mathrm{L}_2(11)$ contains only elements of order 2,3,5,6 and 11 (see \cite{atlas}) and a direct computation shows the following:
\begin{eqnarray}\nonumber
\mathrm{rank}(S_{\alpha}(F_{Kl}))=8&\,\,&\text{if}\,\,\mathrm{ord}(\alpha)=2,\\\nonumber
\mathrm{rank}(S_{\alpha}(F_{Kl}))=12&\,\,&\text{if}\,\,\mathrm{ord}(\alpha)=3,\\\nonumber
\mathrm{rank}(S_{\alpha}(F_{Kl}))=16&\,\,&\text{if}\,\,\mathrm{ord}(\alpha)=5,\\\nonumber
\mathrm{rank}(S_{\alpha}(F_{Kl}))=16&\,\,&\text{if}\,\,\mathrm{ord}(\alpha)=6,\\\nonumber
\mathrm{rank}(S_{\alpha}(F_{Kl}))=20&\,\,&\text{if}\,\,\mathrm{ord}(\alpha)=11.
\end{eqnarray}

\section*{Acknowledgements}
I would like to thank my advisor K. G. O'Grady for his support and for noticing the automorphism $\beta$ in \Ref{sec}{5autom}.\\ I am also grateful to X. Roulleau for pointing out the work contained in \cite{adl} and \cite{kle}, which gave most of the results of \Ref{sec}{5autom}.\\ Moreover I would like to thank A. Rapagnetta for generous advice, M. Sch\"{u}tt for useful discussions and the Referee for helpful comments.


\begin{thebibliography}{}
\bibitem{adl}
A. Adler, On the automorphism group of a certain cubic threefold: Am. Jour. Math. vol. 100 n. 6 (1978) 1275---1280
\bibitem{atlas}
J.H. Conway, R.T. Curtis, S.P. Norton, R.A. Parker and R.A. Wilson, ATLAS of finite groups: Clarendon press, Oxford (1985)
\bibitem{bedon}
A. Beauville and R. Donagi, La variété des droites d'une hypersurface cubique de dimension 4: C.R.Acad.Sc. Paris vol. 301 no. 14 (1985) 703---706 
\bibitem{beau2}
A. Beauville, Some remarks on K\"{a}hler manifolds with $c_1=0$: Prog. Math. vol. 39 (1983) Birkh\"{a}user Boston
\bibitem{con}
J.H. Conway and N.J.A Sloane, Sphere packings, lattices and groups: Grundlehren Math. Wiss. 290 3rd ed. Springer-Verlag (1999)
\bibitem{ero}
V.A. Erokhin, Automorphism groups of 24-dimensional even unimodular lattices (English translation): J. of Sov. Math. vol. 26 (1984) 1876---1879 
\bibitem{kle}
F. Klein, \"{U}ber die Transformationen elfter Ordnung der elliptischen Funktionen: Math. Ann. vol. 15 (1879) 533---555
\bibitem{kon}
S. Kondo, Niemeier lattices, Mathieu groups and finite groups of symplectic automorphisms of $K3$ surfaces: Duke Math. J. vol. 92 no. 3 (1998) 593---603

\bibitem{nik1}
V.V. Nikulin, Finite automorphism groups of K\"{a}hlerian K3 surfaces (Russian): Trudy Moscov Math. Obshch. vol. 38 (1979) 75---137
\bibitem{nik2}
V.V. Nikulin, Integral symmetric bilinear forms and some of their applications (Russian): Izv Akad. Nauk SSSR Ser. Mat. vol. 43 no.1 (1979) 111---177
\bibitem{muk}
S. Mukai, Finite groups of automorphisms of $K3$ surfaces and the Mathieu group: Invent. Math. vol. 94 (1988) 183---221
\bibitem{me1}
G. Mongardi, Symplectic involutions on deformations of $K3^{[2]}$: Centr. Eur. J. Math. vol. 10 no. 4 (2012) 1472---1485 
\bibitem{voi}
C. Voisin, Théorie de Hodge et géométrie algébrique complexe: volume 10 of Cours Spécialisés
[Specialized Courses]. Société Mathématique de France, Paris (2002)
\end{thebibliography}
\end{document}